\documentstyle[11pt,fullpage,twoside,amssym,saks]{article}
\def\Xint#1{\mathchoice
   {\XXint\displaystyle\textstyle{#1}}%
   {\XXint\textstyle\scriptstyle{#1}}%
   {\XXint\scriptstyle\scriptscriptstyle{#1}}%
   {\XXint\scriptscriptstyle\scriptscriptstyle{#1}}%
   \!\int}
\def\XXint#1#2#3{{\setbox0=\hbox{$#1{#2#3}{\int}$}
     \vcenter{\hbox{$#2#3$}}\kern-.5\wd0}}
\def\dashint{\Xint-}

\begin{document}
\markboth{\centerline{V. RYAZANOV, U. SREBRO AND E. YAKUBOV }}
{\centerline{INTEGRAL CONDITIONS IN THE BELTRAMI EQUATIONS THEORY}}
\def\kohta #1 #2\par{\par\noindent\rlap{#1)}\hskip30pt
\hangindent30pt #2\par}
\def\esssup{\operatornamewithlimits{ess\,sup}}
\def\tomes{\mathop{\longrightarrow}\limits^{mes}}
\def\ts{\textstyle}
\def\I{\roman{Im}}
\def\mes{\mbox{\rm mes}}
\def\Rm{{{\Bbb R}^m}}
\def\Rn{{{\Bbb R}^n}}
\def\Rk{{{\Bbb R}^k}}
\def\R3{{{\Bbb R}^3}}
\def\lR{{\overline {{\Bbb R}}}}
\def\lC{{\overline {{\Bbb C}}}}
\def\lz{{\overline {{z}}}}
\def\lD{{\overline {{D}}}}
\def\lR{{\overline {{\Bbb R}}}}
\def\lRn{{\overline {{\Bbb R}^n}}}
\def\lRm{{\overline {{\Bbb R}^m}}}
\def\lBn{{\overline {{\Bbb B}^n}}}
\def\Bn{{{\Bbb B}^n}}
\def\R{{\Bbb R}}
\def\N{{\Bbb N}}
\def\Z{{\Bbb Z}}
\def\C{{\Bbb C}}
\def\B{{\Bbb B}}
\def\Di{{\Bbb D}}
\def\e{{\varepsilon}}
\def\L{{\Lambda}}
\def\f{{\varphi}}
\def\F{{\Phi}}
\def\E{{\eta}}
\def\x{{\chi}}
\def\d{{\delta }}
\def\D{{\Delta }}
\def\c{{\circ }}
\def\tg{{\tilde{\gamma}}}
\def\a{{\alpha }}
\def\b{{\beta }}
\def\p{{\psi }}
\def\P{{\Psi }}
\def\k{{\kappa }}
\def\m{{\mu }}
\def\n{{\nu }}
\def\r{{\rho }}
\def\t{{\tau }}
\def\S{{\Sigma }}
\def\O{{\Omega }}
\def\o{{\omega }}
\def\s{{\sigma }}
\def\v{{\vartheta}}
\def\z{{\zeta }}
\def\l{{\lambda }}
\def\g{{\gamma }}
\def\G{{\Gamma }}
\def\D{{\Delta }}
\let\text=\mbox
\let\Cal=\cal

\def\cc{\setcounter{equation}{0}
\setcounter{figure}{0}\setcounter{table}{0}}

\overfullrule=0pt

\def\eqb{\begin{equation}}
\def\eqe{\end{equation}}
\def\eb{\begin{eqnarray}}
\def\ee{\end{eqnarray}}
\def\ebnn{\begin{eqnarray*}}
\def\eenn{\end{eqnarray*}}
\def\db{\begin{displaystyle}}
\def\de{\end{displaystyle}}
\def\tb{\begin{textstyle}}
\def\te{\end{textstyle}}
\def\exb{\begin{ex}}
\def\exe{\end{ex}}
\def\bth{\begin{theo}}
\def\eth{\end{theo}}
\def\bcor{\begin{corol}}
\def\ecor{\end{corol}}
\def\blem{\begin{lemma}}
\def\elem{\end{lemma}}
\def\brem{\begin{rem}}
\def\erem{\end{rem}}
\def\bpr{\begin{propo}}
\def\epr{\end{propo}}
\title{{\bf INTEGRAL CONDITIONS IN THE THEORY\\OF THE BELTRAMI EQUATIONS}}

\author{{\bf V. Ryazanov, U. Srebro and E. Yakubov}\\ {}\\
DEDICATED TO 85 YEARS OF OLLI LEHTO}
\date{}
\maketitle

\large \abstract It is shown that many recent and new results on the
existence of ACL homeomorphic solutions for the degenerate Beltrami
equations with integral constraints follow from our extension of the
well--known Lehto existence theorem.
\endabstract


\large \cc
\section{Introduction} The classical case was investigated long ago,
see e.g. \cite{Ah}, \cite{Bel}, \cite{Boj} and \cite{LV}. The
existence problem for degenerate Beltrami equations is currently an
active area of research. It has been studied extensively and many
contributions have been made, see e.g. \cite{AIM},
\cite{BGR$_1$}--\cite{BGR$_2$}, \cite{BJ$_1$}--\cite{BJ$_2$},
\cite{Ch}, \cite{Da}, \cite{GMSV$_1$}--\cite{GMSV$_2$},
\cite{IM$_1$}--\cite{IM$_2$}, \cite{Kr}, \cite{Le}, \cite{MM},
\cite{MMV}, \cite{MRSY}, \cite{MS}, \cite{Pe}, \cite{Tu},
\cite{RSY$_1$}--\cite{RSY$_6$} and \cite{Ya}, see also the survey
\cite{SY}. The goal here is to show that our extension of the Lehto
existence theorem has as corollaries the main known existence
theorems as well as a series of more advanced theorems for the
Beltrami equations, see Section 4. The base for these advances is
some lemmas on integral conditions in Sections 2 and 3. Then we show
in Section 5 that the integral conditions found in Section 4 are not
only sufficient but also necessary in the existence theorems for the
Beltrami equations with integral restrictions. Finally, the
corresponding historic comments and final remarks can be found in
Section 6.

Let $D$ be a domain in the complex plane $\C,$ i.e., a connected
open subset of $\C,$ and let $\m: D\to \C$ be a measurable function
with $|\m (z)| < 1$ a.e. The {\bf Beltrami equation} is \eqb
\label{eq1.1} f_{\overline{z}}\, =\, \mu (z)\cdot f_z \eqe where $
f_{\overline{z}} = {\overline {\partial}}f = (f_x+if_y)/2 ,$ $
f_{{z}} = \partial f = (f_x-if_y)/2, $ $ z = x+iy, $ and $f_x$ and
$f_y$ are partial derivatives of $f$ in $x$ and $y,$
correspondingly.

The function $\mu$ is called the {\bf complex coefficient} and \eqb
\label{eq1.1a} K_{\mu}(z)\, =\, \frac{1+|\mu (z)|}{1-|\mu (z)|}\eqe
the {\bf maximal dilatation} or in short the {\bf dilatation} of the
equation (\ref{eq1.1}). The Beltrami equation (\ref{eq1.1}) is said
to be {\bf degenerate} if $ess\, sup\, K_{\m}(z)=\infty .$
\medskip

Use will be made also the {\bf tangential dilatation} with respect
to a point $z_0\in \overline{D}$ which is defined by
\eqb\label{eq1.5P} K^T_{\m}(z,z_0)\ =\
\frac{\left|1-\frac{\overline{z-z_0}}{z-z_0}\mu (z)\right|^2}{1-|\mu
(z)|^2}\ ,  \eqe cf. \cite{An}, \cite{GMSV$_1$}--\cite{GMSV$_2$},
\cite{Le}, \cite{RW}, \cite{RSY$_2$} and \cite{RSY$_3$}. Note the
following precise estimates \eqb \label{eq5.4CC}
\frac{1}{K_{\m}(z)}\ \le\ K^T_{\m}(z,z_0)\ \le\ K_{\m}(z)\ \ \ \ \ \
\ \mbox{a.e.} \eqe Thus, $K^T_{\m}(z,z_0)\ne 0$ and $\infty$ a.e. if
$K_{\m}(z)$ is locally integrable in a domain $D$.
\medskip

Recall that a mapping $f: D\to\C$ is {\bf absolutely continuous on
lines}, abbr. $f\in${\bf ACL}, if, for every closed rectangle $R$ in
$D$ whose sides are parallel to the coordinate axes, $f|R$ is
absolutely continuous on almost all line segments in $R$ which are
parallel to the sides of $R.$ In particular, $f$ is ACL if it
belongs to the Sobolev class $W^{1,1}_{loc},$ see e.g. \cite{Ma}, p.
8. Note that, if $f\in$ ACL, then $f$ has partial derivatives $f_x$
and $f_y$ a.e. Furthermore, every ACL homeomorphism is
differentiable a.e., see e.g. \cite{GL} or \cite{LV}, p. 128, or
\cite{Me} and \cite{Tr}, p. 331. For a sense-preserving ACL
homeomorphism $f: D\to\C ,$ the Jacobian $J_f(z) =
|f_z|^2-|f_{\lz}|^2$ is nonnegative a.e., see \cite{LV}, p. 10. In
this case, the {\bf complex dilatation} of $f$ is the ratio $\mu (z)
= f_{\lz}/f_z,$ and $\,\,\, |\m (z)|\le 1$ a.e., and the {\bf
dilatation} of $f$ is $K_{\m}(z)$ from (\ref{eq1.1a}) and
$K_{\m}(z)\ge 1$ a.e. Here we set by definition $\mu (z) = 0$ and,
correspondingly, $K_{\mu }(z) = 1$ if $f_z = 0.$ The complex
dilatation and the dilatation of $f$ will be denoted by $\m_f$ and
$K_f,$ respectively.
\medskip

 Recall also
that, given a family of paths $\Gamma $ in $\lC ,$ a Borel function
$\rho:\lC \to [0,\infty]$ is called {\bf admissible} for $\Gamma ,$
abbr. $\rho \in adm\, \Gamma ,$ if \eqb \label{eq1.2v}
\int\limits_{\gamma} \rho(z)\, |dz|\ \geq\ 1 \eqe for each
$\gamma\in\Gamma .$ The {\bf modulus} of $\Gamma$ is defined by \eqb
\label{eq1.3v} M(\Gamma) =\inf\limits_{ \rho \in adm\, \Gamma}
\int\limits_{\C} \rho^2(z)\ dxdy\ . \eqe

Given a domain $D$ and two sets $E$ and $F$ in ${\lC}$, $\G (E,F,D)$
denotes the family of all paths $\g:[a,b] \to {\lC}$ which join $E$
and $F$ in $D$, i.e., $\g(a) \in E, \ \g(b) \in F$ and $\g(t) \in D$
for $a<t<b$.  We set $\G(E,F)= \G(E,F,{\lC})$ if $D={\lC}.$ A {\bf
ring domain}, or shortly a {\bf ring} in $\lC$ is a domain $R$ in
$\lC $ whose complement consists of two components. Let $R$ be a
ring in $\lC .$ If $C_1$ and $C_2$ are the components of $\lC
\setminus R,$ we write $R=R(C_1,C_2).$ It is known that $M(\G (C_1,
C_2, R)) = cap \ R(C_1, C_2),$ see e.g. \cite{Ge$_1$}. Note also
that $M(\G (C_1,C_2,R)) = M(\G (C_1,C_2))$, see e.g. Theorem 11.3 in
\cite{Va}. In what follows, we use the notations $B(z_0, r)$ and
$C(z_0, r)$ for the open disk and the circle, respectively, in $\C$
centered at $z_0\in\C$ with the radius $r>0$ and $A(z_0,r_1,r_2)$
for the ring $\{ z\in\C : r_1<|z-z_0|<r_2\}$.

Motivated by the ring definition of quasiconformality in
\cite{Ge$_2$}, we introduced in \cite{RSY$_3$}, cf. also
\cite{RSY$_2$}, the following notion that localizes and extends the
notion of a $Q$--homeomorphism, see e.g. \cite{MRSY}. Let $D$ be a
domain in $\C ,$ $z_0\in D,$ $r_0\le \mbox{dist} (z_0,\partial D)$
and $Q: B(z_0,r _0)\to [0,\infty]$ a measurable function. A
homeomorphism $f:D\to\lC$ is called a {\bf ring $Q$--homeomorphism
at the point} $z_0\in D$ if \eqb \label{eq1.5w} M(\G (fC_1,fC_2))\
\leq \int\limits_{A} Q(z)\cdot \eta^2(|z-z_0|)\ dxdy \eqe for every
ring $A=A(z_0, r_1, r_2),$ $0<r_1<r_2< r_0,$ $C_i=C(z_0, r_i)$,
$i=1,2,$ and for every measurable function $\eta : (r_1,r_2)\to
[0,\infty ]$ such that \eqb \label{eq4.30b}
\int\limits_{r_1}^{r_2}\eta(r)\ dr\ =\ 1\ . \eqe

This notion was first extended to the boundary points in
\cite{RSY$_5$}. More precisely, given a domain $D$ in $\C $ and a
measurable function $Q: D\to [0,\infty]$, we say that a
homeomorphism $f:D\to\lC$ is a {\bf ring $Q-$homeomorphism at a
boundary point} $z_0$ of the domain $D$ if \eqb \label{eq1.555w}
M(\Delta(fC_1,fC_2,fD))\ \leq \int\limits_{A\cap D} Q(z)\cdot
\eta^2(|z-z_0|)\ dxdy \eqe for every ring $A = A(z_0,r_1,r_2)$ and
every continua $C_1$ and $C_2$ in $D$ which belong to the different
components of the complement to the ring $A$ in $\lC$, containing
$z_0$ and $\infty$, correspondingly, and for every measurable
function $\eta : (r_1,r_2)\to [0,\infty ]$ satisfying the condition
(\ref{eq4.30b}).
\medskip

An ACL homeomorphism $f_{\mu}: D \to \C$ is called a {\bf ring
solution} of the Beltrami equation (\ref{eq1.1}) if $f$ satisfies
(\ref{eq1.1}) a.e., $f^{-1} \in W^{1,2}_{loc}(f(D))\ $ and $f$ is a
ring $Q$--homeomorphism at every point $z_0\in D$ with
$Q(z)=K^T_{\m}(z,z_0).$ If in addition $f$ is a ring
$Q$--homeomorphism at every boundary point $z_0\in \partial{D}$ with
$Q(z)=K^T_{\m}(z,z_0)$, then $f$ is called a {\bf strong ring
solution} of (\ref{eq1.1}).
\medskip

The inequality (\ref{eq1.555w}), which strong ring solutions
satisfy, is an useful tool in deriving various, in particular,
boundary properties of such solutions. The condition $f^{-1}\in
W^{1,2}_{loc}$ given in the definition of a ring solution implies
that a.e. point $z$ is a {\bf regular point} for the mapping $f,$
i.e., $f$ is differentiable at $z$ and $J_f(z)\ne 0,$ see \cite{Po}
and Theorem III.6.1 in \cite{LV}. Note that the condition $K_{\m}\in
L^1_{loc}$ is necessary for a homeomorphic ACL solution $f$ of
(\ref{eq1.1}) to have the property $g=f^{-1}\in W^{1,2}_{loc}$
because this property implies that
$$
\int\limits_{C} K_{\m}(z)\ dxdy\ \le\
4\int\limits_{C}\frac{dxdy}{1-|\m (z)|^2}\ =\
4\int\limits_{f(C)}|\partial g|^2\ dudv\ <\ \infty
$$
for every compact set $C\subset D$, see e.g. Lemmas III.2.1, III.3.2
and Theorems III.3.1, III.6.1 in \cite{LV}, cf. also I.C(3) in
\cite{Ah}. Note also that every homeomorphic ACL solution $f$ of the
Beltrami equation with $K_{\m}\in L^1_{loc}$ belongs to the class
$W^{1,1}_{loc}$ as in all our theorems. If in addition $K_{\m}\in
L^p_{loc},$ $p\in [1,\infty ],$ then $f_{\m}\in W^{1,s}_{loc}$ where
$s=2p/(1+p)\in [1,2].$ Indeed, if $f\in$ ACL, then $f$ has partial
derivatives $f_x$ and $f_y$ a.e. and, for a sense-preserving ACL
homeomorphism $f: D\to\C ,$ the Jacobian $J_f(z) =
|f_z|^2-|f_{\lz}|^2$ is nonnegative a.e. and, moreover, \eqb
\label{eq1.6} |\ \overline{\partial}f|\ \leq\ |\ \partial f|\ \leq\
|\
\partial f|\ +\ |\ \overline{\partial}f|\ \leq\ K_{\mu}^{1/2}(z)\cdot
J_f^{1/2}(z)\ \ \ \  \ \ \  a.e.\eqe Recall that if a homeomorphism
$f:D\to\C$ has finite partial derivatives a.e., then \eqb
\label{eq1.7} \int\limits_B\ J_f(z)\ dxdy\ \le\ |f(B)| \eqe for
every Borel set $B\subseteq D$, see e.g. Lemma III.3.3 in \cite{LV}.
Consequently, applying successively the H\"older inequality and the
inequality (\ref{eq1.7}) to (\ref{eq1.6}), we get that \eqb
\label{eq1.8} \Vert\partial f\Vert_s\ \le \ \Vert
K_{\mu}\Vert^{1/2}_p \cdot |f(C)|^{1/2} \eqe where
$\Vert\cdot\Vert_s$ and $\Vert\cdot\Vert_p$ denote the $L^s-$ and
$L^p-$norms in a compact set $C\subset D$, respectively. In the
classical case when $\Vert \m\Vert_{\infty} < 1,$ equivalently, when
$K_{\m}\in L^{\infty},$ every ACL homeomorphic solution $f$ of the
Beltrami equation (\ref{eq1.1}) is in the class $W^{1,2}_{loc}$
together with its inverse mapping $f^{-1}.$ In the case $\Vert
\m\Vert_{\infty} = 1$ and when $K_{\m}\le Q\in \mbox{BMO},$ again
$f^{-1}\in W^{1,2}_{loc}$ and $f$ belongs to $W^{1,s}_{loc}$ for all
$1\le s <2$ but not necessarily to $W^{1,2}_{loc},$ see e.g.
\cite{RSY$_1$}.
\medskip

Olli Lehto  considers in \cite{Le} degenerate Beltrami equations in
the special case where the {\bf singular set} $S_{\m} = \{ z\in \C :
\lim\limits_{\e\to 0}\ \Vert{K_{\m}}\Vert _{L^{\infty}(B(z,\e))} \
=\ \infty \}$ of the complex coefficient $\m$ in (\ref{eq1.1}) is of
measure zero and shows that, if for every $z_0\in\C $, $r_1$ and
$r_2\in (0,\infty)$, $r_2>r_1 ,$ the following integral is positive
and \eqb \label{eq4.30F} \int\limits_{r_1}^{r_2}\
\frac{dr}{r(1+q^T_{z_0}(r))}\ \to \infty \ \ \ \ \ \ \ \ \ \
\mbox{as}\ r_1\to 0\ \mbox{or}\ r_2\to\infty \eqe where
$q^T_{z_0}(r)$ is the average of $K^T_{\m}(z,z_0)$ over $|z-z_0|=r,$
then there exists a homeomorphism $f:\lC\to\lC$ which is ACL in
$\C\setminus S_{\m}$ and satisfies (\ref{eq1.1}) a.e.
\medskip

Our extension and strengthening of Lehto's existence theorem for
ring solutions was first published in the preprint \cite{RSY$_2$},
and then in the journal paper \cite{RSY$_3$}, see also \cite{MRSY}.
The most advanced version for the strong ring solutions is the
following, see \cite{RSY$_6$}:

\bth{} \label{th3.2C}  Let  $D$  be a domain in $\C$
 and let $\mu : D\to\C$ be a measurable function with
 $|\mu (z)| < 1$  a.e.\ and $K_{\m}\in L^1_{loc}(D).$ Suppose that
 \eqb \label{eq4.30C}
 \int\limits_{0}^{\d(z_0)}\frac{dr}{rq^T_{z_0}( r)}\ =\ \infty\ \ \ \ \ \ \ \ \ \ \forall\ z_0\in D \eqe
where $\d(z_0)<dist\, (z_0,\partial D)$ and $q^T_{z_0}(r)$ is the
average of $K^T_{\m}(z,z_0)$ over $|z-z_0|=r.$ Then the Beltrami
equation (\ref{eq1.1}) has a strong ring solution. \eth

Note that the situation where $S_{\m}=D$ is possible here and that
the condition (\ref{eq4.30C}) is a little weaker than the Lehto
condition (\ref{eq4.30F}) because $q^T_{z_0}(r)$ can be arbitrarily
close to $0$, see (\ref{eq5.4CC}). Note also that already in the
work \cite{MS} it was established the existence of
ho\-meo\-mor\-phic solutions to (\ref{eq1.1}) in the class
$f_{\m}\in W^{1,s}_{loc}$, $s=2p/(1+p)$, under the condition
(\ref{eq4.30C}) with $K_{\m}\in L^p_{loc},$ $p>1$, instead of
$K^T_{\m}(z,z_0)$ (the case $p=1$ is covered thanking to a new
convergence theorem in the recent paper \cite{RSY$_5$}, see also
\cite{RSY$_6$}). The Miklyukov-Suvorov result was again discovered
in the paper \cite{Ch} whose author thanks Professor F.W. Gehring
but does not mention \cite{MS}. Perhaps, the work \cite{MS} remains
unknown even for the leading experts in the west because we have
found no reference to this work in the latest monograph in the
Beltrami equation under the discussion of the Lehto condition, see
Theorem 20.9.4 in \cite{AIM}.\bigskip

Theorem 1.14, side by side with lemmas in Sections 2 and 3, is the
the main base for deriving all theorems in Section 4 on existence of
strong ring solutions for Beltrami equations with various integral
conditions.

\cc
\section{On some equivalent integral conditions}

The main existence theorems for the Beltrami equations (\ref{eq1.1})
are based on integral restrictions on the dilatations $K_{\m}(z)$
and $K^T_{\m}(z,z_0)$. Here we establish equivalence of a series of
the corresponding integral conditions.
\medskip

For this goal, we use the following notions of the inverse function
for monotone functions. For every non-decreasing function
$\F:[0,\infty ]\to [0,\infty ] ,$ the {\bf inverse function}
$\F^{-1}:[0,\infty ]\to [0,\infty ]$ can be well defined by setting
\eqb\label{eq5.5CC} \F^{-1}(\tau)\ =\ \inf\limits_{\F(t)\ge \tau}\
t\ . \eqe Here $\inf$ is equal to $\infty$ if the set of
$t\in[0,\infty ]$ such that $\F(t)\ge \t$ is empty. Note that the
function $\F^{-1}$ is non-decreasing, too.

\brem\label{rmk3.333} It is evident immediately by the definition
that \eqb\label{eq5.5CCC} \F^{-1}(\F(t))\ \le\ t\ \ \ \ \ \ \ \
\forall\ t\in[ 0,\infty ] \eqe with the equality in (\ref{eq5.5CCC})
except intervals of constancy of the function $\f(t)$. \erem

Similarly, for every non-increasing function $\f:[0,\infty ]\to
[0,\infty ] ,$ we set \eqb\label{eq5.5C} \f^{-1}(\tau)\ =\
\inf\limits_{\f(t)\le \tau}\ t\ . \eqe Again, here $\inf$ is equal
to $\infty$ if the set of $t\in[0,\infty ]$ such that $\f(t)\le \t$
is empty. Note that the function $\f^{-1}$ is also non-increasing.

\blem{} \label{pr5.AAA} Let $\p:[0,\infty ]\to [0,\infty ]$ be a
sense--reversing ho\-meo\-mor\-phism and $\f:[0,\infty ]\to
[0,\infty ]$ a monotone function. Then \eqb\label{eq5.5Cc}
[\p\circ\f ]^{-1}(\t)\ =\ \f^{-1}\circ\p^{-1}(\t)\ \ \ \ \ \ \
\forall \t\in [0,\infty ] \eqe and \eqb\label{eq5.CCC} [\f\circ\p
]^{-1}(\t)\ \le\ \p^{-1}\circ\f^{-1}(\t)\ \ \ \ \ \ \ \forall \t\in
[0,\infty ] \eqe and, except a countable collection of $\t\in
[0,\infty],$ \eqb\label{eq5.5CC} [\f\circ\p ]^{-1}(\t)\ =\
\p^{-1}\circ\f^{-1}(\t)\ . \eqe The equality (\ref{eq5.5CC}) holds
for all $\t\in [0,\infty]$ iff the function $\f:[0,\infty ]\to
[0,\infty ]$ is strictly monotone. \elem

\brem\label{rmk4.A} If $\p$ is a sense--preserving homeomorphism,
then (\ref{eq5.5Cc}) and (\ref{eq5.5CC}) are obvious for every
monotone function $\f .$ Similar notations and statements also hold
for other segments $[a,b]$, where $a$ and $b\in [-\infty ,+\infty
]$, instead of the segment $[0,\infty]$.\erem

{\it Proof of Lemma \ref{pr5.AAA}.} Let us first prove
(\ref{eq5.5Cc}). If $\varphi$ is non-increasing, then
$$\left[\psi\circ\varphi\right]^{-1}\,(\tau)=\,\inf\limits_{\psi\left(\varphi(t)\right)\geq \tau}t
\,=\,\inf\limits_{\varphi(t)\leq\psi^{-1}(\tau)}t\,=\,\varphi^{-1}\circ\psi^{-1}(\tau)\,.$$
Similarly, if $\varphi$ is non-decreasing, then
$$\left[\psi\circ\varphi\right]^{-1}\,(\tau)=\,
\inf\limits_{\psi\left(\varphi(t)\right)\leq \tau}t
\,=\,\inf\limits_{\varphi(t)\geq\psi^{-1}(\tau)}t\,=\,\varphi^{-1}\circ\psi^{-1}(\tau)\,.$$

Now, let us prove (\ref{eq5.CCC}) and (\ref{eq5.5CC}). If
$\varphi$ is non-increasing, then applying the substitution
$\eta=\psi(t)$ we have
$$\left[\varphi\circ\psi\right]^{-1}\,(\tau)=\,\inf\limits_
{\varphi\left(\psi(t)\right)\geq \tau}t\,=\,
\inf\limits_{\varphi(\eta)\geq \tau}\psi^{-1}(\eta)\,=\,
\psi^{-1}\left(\sup\limits_{\varphi(\eta)\geq
\tau}\eta\right)\,\leq\,$$$$\leq\,\psi^{-1}\left(\inf\limits_{\varphi(\eta)\leq
\tau}\eta\right)\,=\,\psi^{-1}\circ\varphi^{-1}\left(\tau\right)\,,$$
i.e., (\ref{eq5.CCC}) holds for all $\t\in[0,\infty].$ It is evident
that here the strict inequality is possible only for a countable
collection of $\t\in[0,\infty]$ because an interval of constancy of
$\f$ corresponds to every such $\t .$ Hence (\ref{eq5.5CC}) holds
for all $\t\in[0,\infty]$ if and only if $\varphi$ is decreasing.

Similarly, if $\varphi$ is non-decreasing, then
$$\left[\varphi\circ\psi\right]^{-1}\,(\tau)=\,\inf\limits_
{\varphi\left(\psi(t)\right)\leq \tau}t\,=\,
\inf\limits_{\varphi(\eta)\leq \tau}\psi^{-1}(\eta)\,=\,
\psi^{-1}\left(\sup\limits_{\varphi(\eta)\leq
\tau}\eta\right)\,\leq$$$$\leq\,\psi^{-1}\left(\inf\limits_{\varphi(\eta)\geq
\tau}\eta\right)\,=\,\psi^{-1}\circ\varphi^{-1}\left(\tau\right)\
,$$ i.e., (\ref{eq5.CCC}) holds for all $\t\in[0,\infty]$
and again the strict inequality is possible only for a countable collection of $\t\in[0,\infty].$
In the case, (\ref{eq5.5CC}) holds for all $\t\in[0,\infty]$
if and only if $\varphi$ is increasing.

\bcor \label{cor3.C} In particular, if $\f:[0,\infty ]\to [0,\infty
] $ is a monotone function and $\psi = j$ where $j(t)=1/t,$ then
$j^{-1}=j$ and \eqb\label{eq5.5Cd} [j\circ \f ]^{-1}(\t)\ =\
\f^{-1}\circ j(\t)\ \ \ \ \ \ \ \forall \t\in [0,\infty ] \eqe i.e.,
\eqb\label{eq5.5Ce} \f^{-1}(\t)\ =\ \F^{-1}(1/{\t})\ \ \ \ \ \ \
\forall \t\in [0,\infty ] \eqe where $\F=1/{\f},$ \eqb\label{eq5.C}
[\f\circ j ]^{-1}(\t)\ \le\ j\circ\f^{-1}(\t) \ \ \ \ \ \ \ \forall
\t\in [0,\infty ]  \eqe i.e.,  the inverse function of $\f(1/t)$ is
dominated by $1/\f^{-1},$ and except a countable collection of
$\t\in [0,\infty ] $ \eqb\label{eq5.5Cb} [\f\circ j ]^{-1}(\t)\ =\
j\circ\f^{-1}(\t)\ .\eqe $1/\f^{-1}$ is the inverse function of
$\f(1/t)$ if and only if the function $\f $ is strictly monotone.
\ecor

Further, the integral in (\ref{eq333F}) is un\-der\-stood as the
Lebesgue--Stieltjes integral and the integrals in (\ref{eq333Y}) and
(\ref{eq333B})--(\ref{eq333A}) as the ordinary Lebesgue integrals.
In (\ref{eq333Y}) and (\ref{eq333F}) we complete the definition of
integrals by $\infty$ if $\F(t)=\infty ,$ correspondingly,
$H(t)=\infty ,$ for all $t\ge T\in[0,\infty) .$

\bth{} \label{pr4.1aB} Let $\F:[0,\infty ]\to [0,\infty ]$ be a
non-decreasing function and set \eqb\label{eq333E} H(t)\ =\ \log
\F(t)\ .\eqe

Then the equality \eqb\label{eq333Y} \int\limits_{\D}^{\infty}
H'(t)\ \frac{dt}{t}\ =\ \infty   \eqe implies the equality
\eqb\label{eq333F} \int\limits_{\D}^{\infty} \frac{dH(t)}{t}\ =\
\infty  \eqe and (\ref{eq333F}) is equivalent to \eqb\label{eq333B}
\int\limits_{\D}^{\infty}H(t)\ \frac{dt}{t^2}\ =\ \infty \eqe for
some $\D>0,$ and (\ref{eq333B}) is equivalent to every of the
equalities: \eqb\label{eq333C}
\int\limits_{0}^{\d}H\left(\frac{1}{t}\right)\ {dt}\ =\ \infty \eqe
for some $\d>0,$ \eqb\label{eq333D} \int\limits_{\D_*}^{\infty}
\frac{d\E}{H^{-1}(\E)}\ =\ \infty \eqe for some $\D_*>H(+0),$
\eqb\label{eq333A} \int\limits_{\d_*}^{\infty}\ \frac{d\t}{\t
\F^{-1}(\t )}\ =\ \infty \eqe for some $\d_*>\F(+0).$
\medskip

Moreover, (\ref{eq333Y}) is equivalent  to (\ref{eq333F}) and hence
(\ref{eq333Y})--(\ref{eq333A})
 are equivalent each to other  if $\F$ is in addition absolutely continuous.
In particular, all the conditions (\ref{eq333Y})--(\ref{eq333A}) are
equivalent if $\F$ is convex and non--decreasing. \eth

\brem\label{rmk1} It is necessary to give one more explanation. From
the right hand sides in the conditions
(\ref{eq333Y})--(\ref{eq333A}) we have in mind $+\infty$. If
$\Phi(t)=0$ for $t\in[0,t_*]$, then $H(t)=-\infty$ for $t\in[0,t_*]$
and we complete the definition $H'(t)=0$ for $t\in[0,t_*]$. Note,
the conditions (\ref{eq333F}) and (\ref{eq333B}) exclude that $t_*$
belongs to the interval of integrability because in the contrary
case the left hand sides in (\ref{eq333F}) and (\ref{eq333B}) are
either equal to $-\infty$ or indeterminate. Hence we may assume in
(\ref{eq333Y})--(\ref{eq333C}) that $\D>t_0$ where $t_0\colon
=\sup\limits_{\Phi(t)=0}t$, $t_0=0$ if $\F(0)>0$, and $\d<1/t_0$,
correspondingly.   \erem

{\it Proof.} The equality (\ref{eq333Y}) implies (\ref{eq333F})
because except the mentioned special case
$$
\int\limits_{\D}^{T}\ d\Psi(t)\ \geq\ \int\limits_{\D}^{T}
\Psi^{\prime}(t)\ dt\ \ \ \ \ \ \ \ \forall\ T\in(\D,\ \infty)
$$
where
$$
\Psi (t)\ \colon =\ \int\limits_{\D}^t\frac{d\,H(\t)}{\t}\ ,\ \ \
\Psi'(t)\ =\ \frac{H'(t)}{t}\ ,
$$
see e.g. Theorem $7.4$ of Chapter IV in \cite{Sa}, p. 119, and hence
$$
\int\limits_{\D}^{T}\frac{d\,H(t)}{t}\ \geq\ \int\limits_{\D}^{T}
H^{\prime}(t)\ \frac{dt}{t}\ \ \ \ \ \ \ \ \forall\ T\in(\D,\
\infty)
$$

The equality (\ref{eq333F}) is equivalent to (\ref{eq333B}) by
integration by parts, see e.g. Theorem III.14.1  in \cite{Sa}, p.
102. Indeed, again except the mentioned special case, through
integration by parts we have
$$
\int\limits_{\D}^{T}\frac{d\,H(t)}{t}\ -\ \int\limits_{\D}^{T} H(t)\
\frac{dt}{t^2}\ =\ \frac{H(T+0)}{T}\ -\frac{H(\D -0)}{\D}\ \ \ \ \ \
\ \ \forall\ T\in(\D,\ \infty)
$$
and, if
$$
\liminf\limits_{t\to\infty}\ \frac{H(t)}{t}\ < \infty\ ,
$$
then the equivalence of (\ref{eq333F}) and (\ref{eq333B}) is
obvious. If
$$
\lim\limits_{t\to\infty}\ \frac{H(t)}{t}\ = \infty\ ,
$$
then (\ref{eq333B}) obviously holds, $\frac{H(t)}{t}\ge 1$ for
$t>t_0$ and
$$
\int\limits_{t_0}^{T}\frac{d\,H(t)}{t} =
\int\limits_{t_0}^{T}\frac{H(t)}{t}\ \frac{d\,H(t)}{H(t)} \ge \log\
\frac{H(T)}{H(t_0)} = \log\ \frac{H(T)}{T} + \log\ \frac{T}{H(t_0)}
\to \infty
$$
as $T\to\infty ,$ i.e. (\ref{eq333F}) holds, too.
\medskip

Now, (\ref{eq333B}) is equivalent to (\ref{eq333C}) by the change of
variables $t\rightarrow 1/t.$\medskip

Next, (\ref{eq333C}) is equivalent to (\ref{eq333D}) because by the
geometric sense of integrals as areas under graphs of the
corresponding integrands
$$
\int\limits_0^{\d} \Psi(t)\ dt\ =\
\int\limits_{\Psi(\d)}^{\infty}\Psi^{-1}(\eta)\ d\eta\ +\
\d\cdot\Psi(\d)
$$
where $\Psi(t)=H(1/t)$, and because by Corollary \ref{cor3.C} the
inverse function for $H\left(1/t\right)$ coincides with $1/H^{-1}$
at all points except a countable collection.\medskip

Further, set $\psi(\xi)=\log{\xi}.$ Then $H=\psi\circ\Phi$ and by
Lemma \ref{pr5.AAA} and Remark \ref{rmk4.A}
$H^{-1}=\Phi^{-1}\circ\psi^{-1},$ i.e.,
$H^{-1}(\eta)=\Phi^{-1}(e^{\eta}),$ and  by the substitutions
$\tau=e^{\eta},$\,\,$\eta\,=\,\log\ \tau$ we have the equivalence of
(\ref{eq333D}) and (\ref{eq333A}).\medskip

Finally, (\ref{eq333Y}) and (\ref{eq333F}) are equivalent if $\Phi$
is absolutely continuous, see e.g. Theorem IV.7.4 in \cite{Sa} p.
119. \bigskip

\cc
\section{Connection with the Lehto condition}

In this section we establish useful connection of the conditions of
the Lehto type (\ref{eq4.30C}) with one of the integral conditions
from the last section.
\medskip

Recall that a function  $\psi :[0,\infty ]\to [0,\infty ]$ is called
{\bf convex} if $\psi (\lambda t_1 + (1-\lambda) t_2)\le\lambda\psi
(t_1)+ (1-\lambda)\psi (t_2)$ for all $t_1$ and $t_2\in[0,\infty ]$
and $\lambda\in [0,1]$.\medskip

In what follows, $\Di$ denotes the unit disk in the complex plane
$\C$, \eqb\label{eq5.5Cf} \Di\ =\ \{\ z\in\C:\ |z|\ <\ 1\ \}\ .\eqe

\blem{} \label{lem5.5C} Let $Q:\Di\to [0,\infty ]$ be a measurable
function and let $\F:[0,\infty ]\to [0,\infty ]$ be a non-decreasing
convex function. Then \eqb\label{eq3.222} \int\limits_{0}^{1}\
\frac{dr}{rq(r)}\ \ge\ \frac{1}{2}\ \int\limits_{N}^{\infty}\
\frac{d\t}{\t \F^{-1}(\t )} \eqe where $q(r)$ is the average of the
function $Q(z)$ over the circle $|z|=r$ and \eqb\label{eq555.555} N\
=\ \int\limits_{\Di} \F (Q(z))\ dxdy\ .\eqe\elem

{\it Proof.} Note that the result is obvious if $N=\infty$. Hence we
assume further that $N<\infty$. Consequently, we may also assume
that $\F(t)<\infty$ for all $t\in [0,\infty)$ because in the
contrary case $Q\in L^{\infty}(\Di)$ and then the left hand side in
(\ref{eq3.222}) is equal to $\infty$. Moreover, we may assume that
$\F(t)$ is not constant (because in the contrary case
$\F^{-1}(\t)\equiv\infty$ for all $\t>\t_0$ and hence the right hand
side in (\ref{eq3.222}) is equal to 0), $\F(t)$ is (strictly)
increasing, convex and continuous in a segment $[t_*,\infty]$ for
some $t_*\in [0,\infty)$ and \eqb\label{eq5.555Y} \F(t)\ \equiv\
\t_0\ =\ \F(0)\ \ \ \ \ \ \ \ \forall\ t\ \in[0,t_*]\ . \eqe

Next, setting \eqb\label{eq5.555F}H(t)\ \colon =\ \log\ \F(t)\ ,\eqe
we see by Proposition \ref{pr5.AAA} and  Remark \ref{rmk4.A} that
\eqb\label{eq5.555G}H^{-1}(\eta)\ =\ \F^{-1}(e^{\eta})\ ,\ \ \
\F^{-1}(\t)\ =\ H^{-1}(\log\ \t)\ .\eqe Thus, we obtain that
\eqb\label{eq5.555I}q(r) = H^{-1}\left(\log \frac{h(r)}{r^2}\right)
= H^{-1}\left(2\log \frac{1}{r} + \log\ h(r)\right)\ \ \ \ \ \ \ \
\forall\ r\ \in R_*\eqe where $h(r)\ \colon =\ r^2\F(q(r))$ and
$R_*\ =\ \{ r\in(0,1):\ q(r)\ >\ t_*\}$. Then also
\eqb\label{eq5.555K}q(e^{-s})\ =\ H^{-1}\left(2s\ +\ \log\
h(e^{-s})\right)\ \ \ \ \ \ \ \ \forall\ s\ \in S_*\eqe where $ S_*\
=\ \{ s\in(0,\infty):\ q(e^{-s})\ >\ t_*\}$.

Now, by the Jensen inequality \eqb\label{eq5.555L}
\int\limits_0^{\infty} h(e^{-s})\ ds\ =\ \int\limits_0^{1} h(r)\
\frac{dr}{r}\ =\ \int\limits_0^{1} \F(q(r))\ r{dr}\eqe
$$\le\ \int\limits_0^{1}\left( \frac{1}{2\pi}\int\limits_0^{2\pi}
\F(Q(re^{i\vartheta}))\ d\vartheta\right)\ r{dr}\ =\
\frac{N}{2\pi}$$ and then \eqb\label{eq5.555N} |T|\ =\
\int\limits_{T}ds\ \le\ \frac{1}{2}\eqe where $T\  =\ \{\ s\in
(0,\infty):\ \ \ h(e^{-s})\ >\ {N}/{\pi}\}$. Let us show that
\eqb\label{eq5.555O}q(e^{-s})\ \le\ H^{-1}\left(2s\ +\ \log\
\frac{N}{\pi}\right)\ \ \ \ \ \ \ \ \ \ \forall\
s\in(0,\infty)\setminus T_*\eqe where $T_*\ =\ T\cap S_*$. Note that
$(0,\infty)\setminus T_*  = [(0,\infty)\setminus S_*] \cup
[(0,\infty)\setminus T] = [(0,\infty)\setminus S_*] \cup
[S_*\setminus T]$. The inequality (\ref{eq5.555O}) holds for $s\in
S_*\setminus T$ by (\ref{eq5.555K}) because $H^{-1}$ is a
non-decreasing function. Note also that by (\ref{eq5.555Y})
\eqb\label{eq5.K555}  e^{2s}\frac{N}{\pi} = e^{2s} \dashint_{\Di}\
\F(Q(z)) \ dxdy > \F(0) = \t_0 \ \ \ \ \ \ \ \forall\
s\in(0,\infty)\ .\eqe Hence, since the function $\F^{-1}$ is
non-decreasing and $\F^{-1}(\t_0+0)=t_*$, we have by
(\ref{eq5.555G}) that \eqb\label{eq5.L555} t_* <
\F^{-1}\left(\frac{N}{\pi}\ e^{2s}\right) = H^{-1}\left(2s\ +\ \log\
\frac{N}{\pi}\right)\ \ \ \ \ \ \ \ \forall\ s\in(0,\infty)\ .\eqe
Consequently, (\ref{eq5.555O}) holds for $s\in(0,\infty)\setminus
S_*$, too. Thus, (\ref{eq5.555O}) is true.\medskip

Since $H^{-1}$ is non--decreasing, we have by (\ref{eq5.555N}) and
(\ref{eq5.555O}) that \eqb\label{eq5.555P} \int\limits_{0}^{1}\
\frac{dr}{rq(r)}\ =\ \int\limits_{0}^{\infty}\ \frac{ds}{q(e^{-s})}\
\ge\ \int\limits_{(0,\infty)\setminus T_*}\ \frac{ds}{H^{-1}(2s +
\D)}\ \ge\ \eqe
$$
\ge\ \int\limits_{|T_*|}^{\infty}\ \frac{ds}{H^{-1}(2s + \D)}\ \ge\
\int\limits_{\frac{1}{2}}^{\infty}\ \frac{ds}{H^{-1}(2s + \D)}\  =\
\frac{1}{2}\int\limits_{1+\D}^{\infty}\ \frac{d\eta}{H^{-1}(\eta)}
$$
where $ \D=\log N/\pi$. Note that $1+\D = \log\ N + \log\ {e}/{\pi}\
< \log\ N$. Thus, \eqb\label{eq5.555S} \int\limits_{0}^{1}\
\frac{dr}{rq(r)}\ \ge\ \frac{1}{2}\int\limits_{\log N }^{\infty}\
\frac{d\eta}{H^{-1}(\eta)} \eqe and, after the replacement $\eta =
\log\ \t$, we obtain (\ref{eq3.222}).

\bth{} \label{th5.555} Let $Q:\Di\to [0,\infty ]$ be a measurable
function such that \eqb\label{eq5.555} \int\limits_{\Di} \F (Q(z))\
dxdy\  <\ \infty\eqe where $\F:[0,\infty ]\to [0,\infty ]$ is a
non-decreasing convex function such that \eqb\label{eq3.333a}
\int\limits_{\d_0}^{\infty}\ \frac{d\t}{\t \F^{-1}(\t )}\ =\ \infty
\eqe for some $\d_0\ >\ \t_0\ \colon =\ \F(0).$ Then
\eqb\label{eq3.333A} \int\limits_{0}^{1}\ \frac{dr}{rq(r)}\ =\
\infty \eqe where $q(r)$ is the average of the function $Q(z)$ over
the circle $|z|=r$. \eth

\brem\label{rmk4.7www}  Note that (\ref{eq3.333a}) implies that
\eqb\label{eq3.a333} \int\limits_{\d}^{\infty}\ \frac{d\t}{\t
\F^{-1}(\t )}\ =\ \infty \eqe for every $\d\ \in\ [0,\infty)$ but
(\ref{eq3.a333}) for some $\d\in[0,\infty)$, generally speaking,
does not imply (\ref{eq3.333a}). Indeed, for $\d\in [0,\d_0),$
(\ref{eq3.333a}) evidently implies (\ref{eq3.a333})  and, for
$\d\in(\d_0,\infty)$, we have that \eqb\label{eq3.e333} 0\ <\
\int\limits_{\d_0}^{\d}\ \frac{d\t}{\t \F^{-1}(\t )}\ \le\
\frac{1}{\F^{-1}(\d_0)}\ \log\ \frac{\d}{\d_0}\ <\ \infty \eqe
because $\F^{-1}$ is non-decreasing and $\F^{-1}(\d_0)>0$. Moreover,
by the definition of the inverse function $\F^{-1}(\t)\equiv 0$ for
all $\t \in [0,\t_0],$ $\t_0=\F(0)$, and hence (\ref{eq3.a333}) for
$\d\in[0,\t_0),$ generally speaking, does not imply
(\ref{eq3.333a}). If $\t_0 > 0$, then \eqb\label{eq3.c333}
\int\limits_{\d}^{\t_0}\ \frac{d\t}{\t \F^{-1}(\t )}\ =\ \infty\ \ \
\ \ \ \ \ \ \ \forall\ \d\ \in\ [0,\t_0) \eqe However,
(\ref{eq3.c333}) gives no information on the function $Q(z)$ itself
and, consequently, (\ref{eq3.a333}) for $\d < \F(0)$ cannot imply
(\ref{eq3.333A}) at all. \erem

By  (\ref{eq3.a333}) the proof of Theorem \ref{th5.555} is reduced
to Lemma \ref{lem5.5C}. Combining Theorems \ref{pr4.1aB} and
\ref{th5.555} we also obtain the following conclusion.

\bcor \label{cor555} If $\F:[0,\infty ]\to [0,\infty ]$ is a
non-decreasing convex function and $Q$ satisfies the condition
(\ref{eq5.555}), then every of the conditions
(\ref{eq333Y})--(\ref{eq333A}) implies (\ref{eq3.333A}). \ecor

\cc
\section{Sufficient conditions for solvability}

The following existence theorem is obtained immediately  from
Theorems \ref{th3.2C} and \ref{th5.555}.

\bth{} \label{th4.6c} Let $\mu : D\to\C$ be a measurable function
with $|\mu (z)| < 1$  a.e. and $K_{\m}\in L^1_{loc}.$ Suppose that
every point $z_0\in D$ has a neighborhood $U_{z_0}$ where \eqb
\label{eq5.26V} \int\limits_{U_{z_0}}\ \F_{z_0}(K^T_{\m}(z,z_0))\
dxdy\ <\ \infty \eqe for a non-decreasing convex function $\F_{z_0}
:[0,\infty )\to [0,\infty ]$ such that \eqb\label{eq5.26La}
\int\limits_{\D(z_0)}^{\infty}\frac{d\t}{\t\F^{-1}_{z_0}(\t)}\ \ =\
\infty \eqe for some $\D(z_0)>\F_{z_0}(0).$ Then the Beltrami
equation (\ref{eq1.1}) has a strong ring solution. \eth

{\it Proof.} Let the closure of a disk $B(z_0,\r)$ belong to the
neighborhood $U_{z_0}$. Then we obtain by Theorem \ref{th5.555}
applied to  $Q(\z)=K^T_{\m}(z_0+\r\, \z, z_0),$ $\z\in\Di$, and
$\F(t)=\F_{z_0}(t)$ that \eqb\label{eq3.222a} \int\limits_{0}^{\r}\
\frac{dr}{rq^T_{z_0}(r)}\ =\ \infty \eqe where $q^T_{z_0}(r)$ is the
mean value of $K^T_{\m}(z,z_0)$ over the circle $|z-z_0|=r.$ Thus,
 we have the desired conclusion by Theorem \ref{th3.2C}.

\brem\label{rmk4.13A} Note that the additional condition
$\D(z_0)>\F_{z_0}(0)$ is essential, see also Remark \ref{rmk4.7www}.
In fact, it is important only degree of convergence
$\F^{-1}_{z_0}(\t)\to\infty$ as $\t\to \infty$ or, the same, degree
of convergence $\F_{z_0}(t)\to \infty$ as $t\to\infty .$\erem

\bcor \label{cor4.6C} Let $\mu : D\to\C$ be a measurable function
with $|\mu (z)| < 1$  a.e. and $K_{\m}\in L^1_{loc}.$ Suppose that
\eqb \label{eq5.26VVV} \int\limits_{D}\ \F(K_{\m}(z))\ dxdy\ <\
\infty \eqe for a non-decreasing convex function $\F :[0,\infty ]\to
[0,\infty ]$ such that \eqb\label{eq5.26Laa}
\int\limits_{\D}^{\infty}\frac{d\t}{\t\F^{-1}(\t)}\ \ =\ \infty \eqe
for some $\D>\F(0).$ Then the Beltrami equation (\ref{eq1.1}) has a
strong ring solution. \ecor

\brem\label{rmk4.7s} Applying transformations $\a\cdot\F + \b$ with
$\a>0$ and $\b\in\R$, we may assume without loss of generality that
$\F(t)=\F(1)=1$ for all $t\in [0,1]$ and, thus, $\F(0)=\F(1)=1$ in
Theorem \ref{th4.6c} and its corollaries further. \erem

Many other criteria of the existence of strong ring solutions for
the Beltrami equation (\ref{eq1.1}) formulated below follow from
Theorems \ref{pr4.1aB} and  \ref{th4.6c}.

\bcor \label{cor333} Let $\mu : D\to\C$ be a measurable function
with $|\mu (z)| < 1$  a.e. and $K_{\m}\in L^1_{loc}$. If the
condition (\ref{eq5.26V}) holds at every point $z_0\in D$ with a
non-decreasing convex function $\F_{z_0} :[0,\infty )\to [0,\infty
)$ such that \eqb\label{eq300} \int\limits_{\D(z_0)}^{\infty} \log
\F_{z_0}(t)\ \frac{dt}{t^2}\ =\ \infty \eqe for some $\D(z_0)>0,$
then (\ref{eq1.1}) has a strong ring solution.\ecor

\bcor \label{cor300} Let $\mu : D\to\C$ be a measurable function
with $|\mu (z)| < 1$  a.e. and $K_{\m}\in L^1_{loc}$. If the
condition (\ref{eq5.26V}) holds at every point $z_0\in D$ for a
continuous non--decreasing convex function $\F_{z_0} :[0,\infty )\to
[0,\infty )$ such that \eqb\label{eq5.6a}
\int\limits_{\D(z_0)}^{\infty}(\log \F_{z_0}(t))'\ \frac{dt}{t}\ =\
\infty \eqe for some $\D(z_0)>0,$ then (\ref{eq1.1}) has a strong
ring solution.\ecor

\bcor \label{cor4.6d} Let $\mu : D\to\C$ be a measurable function
with $|\mu (z)| < 1$  a.e. and $K_{\m}\in L^1_{loc}.$ If the
condition (\ref{eq5.26V}) holds at every point $z_0\in D$ for
$\F_{z_0}=\exp{ H_{z_0}}$ where $H_{z_0}$ is non--constant,
non--decreasing and convex, then (\ref{eq1.1}) has a strong ring
solution.\ecor

\bcor \label{cor4.6e} Let $\mu : D\to\C$ be a measurable function
with $|\mu (z)| < 1$  a.e. and $K_{\m}\in L^1_{loc}$. If the
condition (\ref{eq5.26V}) holds at every point $z_0\in D$ for
$\F_{z_0}=\exp H_{z_0}$ with a twice continuously differentiable
increasing function $H_{z_0}$ such that \eqb\label{eq5.6d}
H''_{z_0}(t))\ \ge\  0\ \ \ \ \ \ \ \ \ \forall\ t\ \ge\ t(z_0)\
\in\ [0,\infty)\ , \eqe then (\ref{eq1.1}) has a strong ring
solution. \ecor

\bth{} \label{th4.5A} Let $\mu : D\to\C$ be a measurable function
with $|\mu (z)| < 1$  a.e. and $K_{\m}\in L^1_{loc}$ such that \eqb
\label{eq5.27V} \int\limits_D\ \F(K_{\m}(z))\ dxdy\ <\ \infty \eqe
where $\F :[0,\infty )\to [0,\infty ]$ is non-decreasing and convex
such that \eqb\label{eq301} \int\limits_{\D}^{\infty} \log \F(t)\
\frac{dt}{t^2}\ =\ \infty \eqe for some $\D>0.$ Then the Beltrami
equation (\ref{eq1.1}) has a strong ring solution. \eth

\bcor \label{cor301} Let $\mu : D\to\C$ be a measurable function
with $|\mu (z)| < 1$  a.e. and $K_{\m}\in L^1_{loc}$. If the
condition (\ref{eq5.27V}) holds  with a non--decreasing convex
function $\F :[0,\infty )\to [0,\infty )$ such that
\eqb\label{eq5.27w} \int\limits_{t_0}^{\infty}(\log \F(t))'\
\frac{dt}{t}\ =\ \infty \eqe for some $t_0>0,$ then (\ref{eq1.1})
has a strong ring solution.\ecor

\bcor \label{cor4.6v} Let $\mu : D\to\C$ be a measurable function
with $|\mu (z)| < 1$  a.e. If the condition (\ref{eq5.27V}) holds
for $\F=e^H$ where $H$ is non--constant, non--decreasing and convex,
then (\ref{eq1.1}) has a strong ring solution.\ecor

\bcor \label{cor4.6w} Let $\mu : D\to\C$ be a measurable function
with $|\mu (z)| < 1$  a.e. If the condition (\ref{eq5.27V}) holds
for $\F=e^H$ where $H$ is twice continuously differentiable,
increasing and \eqb\label{eq5.6d} H''(t)\ \ge\ 0\ \ \ \ \ \ \ \
\forall\ t\ \ge\ t_0\ \in\ [1,\infty)\ , \eqe then (\ref{eq1.1}) has
a strong ring solution. \ecor

Note that among twice continuously  differentiable functions, the
condition (\ref{eq5.6d}) is equivalent to the convexity of $H(t),$
$t\ge t_0$, cf. Corollary \ref{cor4.6v}. Of course, the convexity of
$H(t)$ implies the convexity of $\F(t)=e^{H(t)},$ $t\ge t_0,$
because the function $\exp x$ is convex. However, in general, the
convexity of $\F$ does not imply the convexity of $H(t)=\log \F(t)$
and it is known that the convexity of $\F(t)$ in (\ref{eq5.27V}) is
not sufficient for the existence of ACL homeomorphic solutions of
the Beltrami equation. There exist examples of the complex
coefficients $\m$ such that $K_{\m}\in L^p$ with an arbitrarily
large $p\ge 1$ for which the Beltrami equation (\ref{eq1.1}) has no
ACL homeomorphic solutions, see e.g. \cite{RSY$_1$}.
\medskip

\brem\label{rmk111} Theorem \ref{th3.2C} is extended by us to the
case where $\infty\in D\subset\lC $ in the standard way by replacing
(\ref{eq4.30C}) to the following condition at $\infty $ \eqb
\label{eq4.30R}
 \int\limits_{\d}^{\infty}\frac{dr}{rq( r)}\ =\ \infty \eqe
 where $\d > 0$ and  \eqb
\label{eq4.30P} q(r)\ =\ \frac{1}{2\pi}\
\int\limits_{0}\limits^{2\pi}\
\frac{|1-e^{-2i\v}\m(re^{i\v})|^2}{1-|\m(re^{i\v})|^2}\ \ d\v\ .
\eqe In this case, there exists a homeomorphic $W^{1,1}_{loc}$
solution $f$ of (\ref{eq1.1}) in $D$ with $f(\infty)=\infty $ and
$f^{-1}\in W^{1,2}_{loc}.$ Here $f\in W^{1,1}_{loc}$ in $D$ means
that $f\in W^{1,1}_{loc}$ in $D\setminus\{\infty\}$ and that
$f^*(z)=1/\overline{f(1/\lz)}$ belongs to $ W^{1,1}$ in a
neighborhood of $0.$ The statement $f^{-1}\in W^{1,2}_{loc}$ has a
similar meaning.
\medskip

Similarly, the integral condition (\ref{eq5.26V}) is replaced at
$\infty$  by the following condition  \eqb \label{eq5.26q}
\int\limits_{|z|> \d}\ \F_{\infty}(K^T_{\m}(z,\infty))\
\frac{dxdy}{|z|^4}\ <\ \infty \eqe where $\d > 0,$ $\F_{\infty}$
satisfies the conditions of either Theorem \ref{th4.6c} or
equivalent conditions from Theorem \ref{pr4.1aB} and
\eqb\label{eq1.5L} K^T_{\m}(z,\infty)\ =\
\frac{\left|1-\frac{\overline{z}}{z}\mu (z)\right|^2}{1-|\mu
(z)|^2}\ .  \eqe

We may assume in all the above theorems that the functions
$\F_{z_0}(t)$ and $\F(t)$ are not convex on the whole segments
$[0,\infty]$ and $[1,\infty]$, respectively, but only on a segment
$[T,\infty]$ for some $T\in(1,\infty)$. Indeed, every non-decreasing
function $\F:[1,\infty]\to[0,\infty]$ which is convex on a segment
$[T,\infty]$, $T\in(0,\infty)$, can be replaced by a non-decreasing
convex function $\F_T:[0,\infty]\to[0,\infty]$ in the following way.
We set $\F_T(t)\equiv 0$ for all $t\in [0,T]$, $\F(t)=\f(t)$,
$t\in[T,T_*]$, and $\F_T\equiv \F(t)$, $t\in[T_*,\infty]$, where
$\t=\f(t)$ is the line passing through the point $(0,T)$ and
touching upon the graph of the function $\t=\F(t)$ at a point
$(T_*,\F(T_*))$, $T_*\ge T$. For such a function we have by the
construction that $\F_T(t)\le \F(t)$ for all $t\in[1,\infty]$ and
$\F_T(t)=\F(t)$ for all $t\ge T_*$. \erem

\cc
\section{Necessary conditions for solvability}

The main idea for the proof of the following statement under smooth
increasing functions $\F$ with the additional condition that
$t(\log\F)'\ge 1$ is due to Iwaniec and Martin, see Theorem 3.1 in
\cite{IM$_2$}, cf. also Theorem 11.2.1 in \cite{IM$_1$} and Theorem
20.3.1 in \cite{AIM}. We obtain the same conclusion in Theorem
\ref{th55} and Lemma \ref{th555} further without these smooth
conditions. Moreover, by Theorem \ref{pr4.1aB} the same conclusion
concerns to all conditions (\ref{eq333F})--(\ref{eq333A}).

\bth{} \label{th55} Let $\F : [0,\infty]\to[0,\infty]$ be such a
non-decreasing convex function that, for every measurable function
$\m : \Di\to\Di$ satisfying the condition \eqb \label{eq5.27VVV}
\int\limits_{\Di}\ \F(K_{\m}(z))\ dxdy\ <\ \infty\ , \eqe   the
Beltrami equation (\ref{eq1.1}) has a homeomorphic ACL solution.
Then there is $\delta
>0$ such that \eqb\label{eq29.29A}
\int\limits_{\delta}^{\infty}\log{\Phi(t)}\ \frac{dt}{t^2}\ =\
\infty\ . \eqe \eth

It is evident that the function $\Phi(t)$ in Theorem \ref{th55} is
not constant on $[0,\infty)$ because in the contrary case we would
have no real restrictions for $K_{\m}$ from (\ref{eq5.27VVV}) except
$\Phi(t)\equiv\infty$ when the class of such $\mu$ is empty.
Moreover, by the well--known criterion of convexity, see e.g.
Proposition 5 in I.4.3 of \cite{Bou}, the inclination
$[\Phi(t)-\Phi(0)]/t$ is nondecreasing. Hence the proof of Theorem
\ref{th55} is reduced to the following statement.

\blem{} \label{th555} Let a function $\F : [0,\infty]\to[0,\infty]$
be non-decreasing and \eqb\label{eq3.VVV} \F(t)\ \ge\ {C}\cdot{t}\ \
\ \ \ \ \ \ \forall\ t\ \ge\ T \eqe for some $C>0$ and
$T\in(1,\infty)$. If the Beltrami equations (\ref{eq1.1}) have ACL
homeomorphic solutions for all measurable functions $\m : \Di\to\Di$
satisfying the condition (\ref{eq5.27VVV}), then (\ref{eq29.29A})
holds for some $\delta >0$. \elem

\brem\label{rmk6.333} Note that the Iwaniec--Martin condition
$t(\log\F)'\ge 1$ implies the condition (\ref{eq3.VVV}) with
$C=\Phi(T)/T$. Note also that if we take further in the construction
of Lemma \ref{lem6.A} $\b_{n+1}=\e_n\to 0$ as $n\to\infty$,
$\a_{n+1}=b_n e^{b_n\gamma_n}-\e_n\g_n$ and $\g^*_{n+1}=b_n
e^{b_n\gamma_n}/\e_n,$ then $\F(\g^*_{n+1})/\g^*_{n+1}\le 2\e_n$ and
we obtain examples of absolutely continuous increasing functions
$\F$ with $\F(t)\to\infty$ as $t\to\infty$ satisfying the condition
(\ref{eq29.29A}) and simultaneously \eqb\label{eq300DDD}
\liminf\limits_{t\to\infty} \ \frac{\F(t)}{ t}\ =\ 0\ .\eqe Thus,
conditions of the type (\ref{eq3.VVV}) are independent on the
conditions (\ref{eq333Y})--(\ref{eq333A}). \erem

{\it Proof of Lemma \ref{th555}}. Let us assume that the condition
(\ref{eq29.29A}) does not hold for any $\delta > 0$. Set
$t_0=\sup\limits_{\Phi(t)=0}t$, $t_0=0$ if $\F(t)>0$ for all
$t\in[0,\infty]$.  Then for all $\d>t_0$ \eqb\label{eq29.29}
\int\limits_{\delta}^{\infty}\log{\Phi(t)}\ \frac{dt}{t^2}\ <\
\infty\ . \eqe

With no loss of generality, applying the linear transformation $\a\F
+ \b$ with $\a = 1/C$ and $\b = T$,  we may assume by
(\ref{eq3.VVV}) that \eqb\label{eqKKK3} \F(t)\ \ge\ t\ \ \ \ \ \ \ \
\ \ \ \ \forall\ t\in[0,\infty)\ .\eqe Of course, we may also assume
that $\F(t)=t$ for all $t\in[0,1)$ because the values of $\F$ in
$[0,1)$ give no information on $K_{\mu}$ in (\ref{eq5.27VVV}).
Finally, by (\ref{eq29.29}) we have that $\Phi(t)<\infty$ for every
$t\in[0,\infty)$.

Now, note that the function $\P(t)\colon = t\F(t)$ is strictly
increasing, $\P(1)=\F(1)$ and $\P(t)\to\infty$ as $t\to\infty$,
$\P(t)<\infty$ for every $t\in[0,\infty)$. Hence the functional
equation \eqb\label{eqLLL3} \P(K(r))\ =\
\left(\frac{\g}{r}\right)^2\ \ \ \ \ \ \ \ \ \ \ \ \ \forall\ r\
\in\ (0,1]\ ,\eqe where $\g=\F^{1/2}(1)\ge 1$, is well solvable with
$K(1)=1$ and a continuous non-increasing function $K : (0,1]\to
[1,\infty)$ such that $K(r)\to\infty$ as $r\to 0.$ Taking the
logarithm in (\ref{eqLLL3}), we have that$$ 2\log\ r\ +\ \log\ K(r)\
+\ \log\ \F(K(r))\ =\ 2\ log\ \g $$ and by (\ref{eqKKK3}) we obtain
that $$ \log\ r\ +\ \log\ K(r)\ \le\ log\ \g\ , $$ i.e.,
\eqb\label{eqMMM3} K(r)\ \le\ \frac{\g}{r}\ . \eqe Then by
(\ref{eqLLL3})
$$
\F(K(r))\ \ge\ \frac{\g}{r}
$$
and hence by (\ref{eq5.5CCC})
$$
K(r)\ \ge\ \F^{-1}\left(\frac{\g}{r}\right)\ .
$$
Thus,
$$
I(t)\ \colon =\ \int\limits_0^t\ \frac{dr}{rK(r)}\ \le\
\int\limits_0^t\ \frac{dr}{r\F^{-1}\left(\frac{\g}{r}\right)}\ =\
\int\limits_{\frac{\g}{t}}^{\infty}\ \frac{d\t}{\t\F^{-1}(\t)}, \ \
\ \ \ \ \ t\in(0,1]\ ,
$$
where $\g/t \ge \g \ge 1 > \F(+0)=0$. Hence by the condition
(\ref{eq29.29}) and Proposition \ref{pr4.1aB} \eqb\label{eqNNN3}
I(t)\ \le\ I(1)\ =\ \int\limits_0^1\ \frac{dr}{rK(r)}\ <\ \infty\ .
\eqe

Next, consider the mapping
$$
f(z)\ =\ \frac{z}{|z|}\ \r(|z|)
$$
where $\r(t)= e^{I(t)}$. Note that $f\in C^1(\Di\setminus\{ 0\})$
and hence $f$ is locally quasiconformal in the punctured unit disk
$\Di\setminus\{ 0\}$ by the continuity of the function $K(r),$
$r\in(0,1)$, see also (\ref{eqMMM3}). Let us calculate its complex
dilatation. Set $z=re^{i\vartheta}$. Then
$$
\frac{\partial f}{\partial r}\ =\ \frac{\partial f}{\partial z}
\cdot \frac{\partial z}{\partial r}\ +\ \frac{\partial f}{\partial
\overline {z}} \cdot \frac{\partial \overline {z}}{\partial r}\ =\
e^{i\vartheta}\cdot\frac{\partial f}{\partial z} \ +\
e^{-i\vartheta}\cdot\frac{\partial f}{\partial \overline {z}}
$$
and
$$ \frac{\partial f}{\partial \vartheta}\ =\ \frac{\partial
f}{\partial z} \cdot \frac{\partial z}{\partial \vartheta}\ +\
\frac{\partial f}{\partial \overline {z}} \cdot \frac{\partial
\overline {z}}{\partial \vartheta}\ =\
ire^{i\vartheta}\cdot\frac{\partial f}{\partial z} \ -\
ire^{-i\vartheta}\cdot\frac{\partial f}{\partial \overline {z}}\ .
$$
In other words, \eqb\label{eqFFF3} \frac{\partial f}{\partial z}\ =\
\frac{e^{-i\vartheta}}{2}\left(\frac{\partial f}{\partial r} \ +\
\frac{1}{ir}\cdot\frac{\partial f}{\partial \vartheta}\right) \eqe
and \eqb\label{eqEEE3} \frac{\partial f}{\partial \overline {z}}\ =\
\frac{e^{i\vartheta}}{2}\left(\frac{\partial f}{\partial r} \ -\
\frac{1}{ir}\cdot\frac{\partial f}{\partial \vartheta}\right)\eqe

Thus, we have that

$$ \frac{\partial f}{\partial z}\ =\
\frac{1}{2}\left(\frac{\r (r)}{rK(r)} \ +\ \frac{\r(r)}{r}\right)\
=\ \frac{\r(r)}{2r}\ \cdot \frac{1+K(r)}{K(r)}
$$
and
$$ \frac{\partial f}{\partial \overline {z}}\ =\
\frac{e^{2i\vartheta}}{2}\left(\frac{\r (r)}{rK(r)} \ -\
\frac{\r(r)}{r}\right)\ =\ e^{2i\vartheta}\cdot\frac{\r(r)}{2r}\
\cdot \frac{1-K(r)}{K(r)}
$$
i.e.
$$ \m(z)\ =\
e^{2i\vartheta}\cdot\frac{1-K(r)}{1+K(r)}\ =\ -\ \frac{z}{\overline
{z}}\cdot\frac{K(|z|)-1}{K(|z|)+1}\ .
$$
Consequently, \eqb\label{eqOOO3} K_{\m}(z)\ =\ K(|z|)\eqe and by
(\ref{eqLLL3})
$$
 \int\limits_{\Di}\ \F(K_{\m}(z))\ dxdy\ =\ 2\pi\ \int\limits_{0}^1\ \F(K(r))\ r\ dr\
 \le\ 2\pi\g^2 I(1)\ <\ \infty\ .
$$
However,
$$
\lim\limits_{z\to 0}\ |f(z)|\ =\ \lim\limits_{t\to 0}\ \r(t)\ =\
e^{I(0)}\ =\ 1\ ,
$$
i.e. $f$ maps the punctured disk $\Di\setminus\{ 0\}$ onto the ring
$1<|\z|< R=e^{I(1)}$.\medskip

Let us assume that there is a homeomorphic ACL solution $g$ of the
Beltrami equation (\ref{eq1.1}) with the given $\m$.  By the Riemann
theorem without loss of generality we may assume that $g(0)=0$ and
$g(\Di)=\Di$. Since $f$ as well as $g$ are locally quasiconformal in
the punctured disk $\Di\setminus\{ 0\}$, then by the uniqueness
theorem for the quasiconformal mappings $f=h\circ g$ in
$\Di\setminus\{ 0\}$ where $h$ is a conformal mapping in
$\Di\setminus\{ 0\}$. However, isolated singularities are removable
for conformal mappings. Hence $h$ can be extended  by continuity to
$0$ and, consequently, $f$ should be so. Thus, the obtained
contradiction disproves the assumption (\ref{eq29.29}).

\brem\label{rmk555} Thus, Theorems \ref{th55} and \ref{pr4.1aB} show
that every of the con\-di\-tions (\ref{eq333Y})--(\ref{eq333A}) in
the existence theorems  to the Beltrami equations (\ref{eq1.1}) with
the integral constraints (\ref{eq5.27VVV}) for non--decreasing
convex functions $\F$ are not only sufficient but also necessary.
\erem

\cc
\section{Historic comments and final remarks} To compare our results
with erlier results of other authors we give a short survey.\bigskip

The first investigation of the existence problem for degenerate
Beltrami equa\-tions with integral constraints (\ref{eq5.27V}) as
in Theorem \ref{th4.5A} has been made by Pesin \cite{Pe} who
studied the special case where $\F(t)=e^{t^{\a}}-1$ with $\a
> 1.$ Basically, Corollary \ref{cor4.6w} is due to Kruglikov
\cite{Kr}. David \cite{Da} considered the existence problem with
measure constraints \eqb \label{eq5.26O} |\{ z\in D:\ K_{\m}(z)\ >\
t\}|\ \le\ \f (t)\ \ \ \ \ \ \ \ \forall\ t\in [1,\infty)\eqe with
special $\f(t)$ of the form $a\cdot e^{-b t}$  and Tukia \cite{Tu}
with the corresponding constraints in terms of the spherical area.
Note that under the integral constraints (\ref{eq5.27V}) of the
exponential type $\F(t)=\a e^{\b t}$, $\a > 0$, the conditions of
David and Tukia hold. Thus,  the latter results sthrengthen  the
Pesin result.\medskip

By the well known John-Nirenberg lemma for the function of the class
BMO (bounded mean oscillation) the David conditions are equivalent
to the cor\-res\-pon\-ding integral conditions of the exponential
type, see e.g. \cite{RSY$_1$}. More advanced results in terms of FMO
(finite mean oscillation by Ignat'ev-Ryazanov) can be found in
\cite{RSY$_4$}, \cite{RSY$_6$} and \cite{MRSY}.
\medskip

The next step has been made by Brakalova and Jenkins \cite{BJ$_1$}
who proved the existence of ACL homeomorphic solutions for the case
of the integral constraints (\ref{eq5.26V}) as in Theorem
\ref{th4.6c} with $K_{\m}(z)$ instead of $K_{\m}^T(z,z_0)$ and with
\eqb\label{eq55a} \F_{z_0}(t)\ \equiv\ \F(t)\ =\ \exp\
\left(\frac{\frac{t+1}{2}}{1+\log\frac{t+1}{2}}\right)\ . \eqe Note
that, in the case \cite{BJ$_1$}, the condition (\ref{eq333Y}) in
Theorem \ref{pr4.1aB}, see also Corollary \ref{cor300}, can be easy
verified by the calculations \eqb\label{eq55b}
 (\log\F(t))'\ =\
\frac{1}{2}\frac{\log\frac{t+1}{2}}{(1+\log\frac{t+1}{2})^2}\ \sim\
\frac{1}{2}\frac{1}{\log t}\ \ \ \ \ \ \mbox{as}\ t\to\ \infty\
.\eqe Moreover, it is easy to verify that $\F''(t)\ge 0$ for all
$t\ge T$ under large enough $T\in(1,\infty)$ and, thus, $\F$ is
convex on the segment $[T,\infty]$, see e.g. \cite{Bou} and Remark
\ref{rmk111}.

Later on, Iwaniec and Martin have proved the existence of solutions
in the Orlicz--Sobolev classes for the case where \eqb\label{eq555}
\F_{z_0}(t)\ \equiv\ \F(t)\ =\ \exp\ \left(\frac{p t}{1+\log
{t}}\right)\eqe for some $p>0,$ see e.g.
\cite{IM$_1$}--\cite{IM$_2$}, for which \eqb\label{eq555b}
 (\log\F(t))'\ =\
\frac{p\log\ t}{(1+\log{t})^2}\ \sim\ \frac{p}{\log t}\ \ \ \ \ \
\mbox{as}\ t\to\ \infty\ ,\eqe cf. Corollary \ref{cor301}. Note that
in the both cases (\ref{eq55a}) and (\ref{eq555}) \eqb\label{eq555c}
\F(t)\ \ge\ t^{\l}\ \ \ \ \ \ \ \ \forall\ t\ \ge\ t_{\l}\ \in\
[1,\infty)\ . \eqe It is remarkable that in the case
\eqb\label{eq555.bmo} \F_{z_0}(t)\ \equiv\ \F(t)\ =\ \exp\ {p t}\eqe
it was established uniqueness and factorization theorems for
solutions of the Beltrami equations of the Stoilow type, see e.g.
\cite{AIM} and \cite{Da}.
\medskip

Corollary \ref{cor333} is due to Gutlyanskii, Martio, Sugawa and
Vuorinen in \cite{GMSV$_1$} and \cite{GMSV$_2$} where they have
established the existence of ACL homeomorphic solutions of
(\ref{eq1.1}) in $W^{1,s}_{loc}$, $s=2p/(1+p)$, under $K_{\m}\in
L^p_{loc}$ with $p>1$  for \eqb\label{eq56a} \F_{z_0}(t)\ \equiv\
\F(t)\ \colon =\ \exp H(t) \eqe with $H(t)$ being a continuous
non--decreasing function such that $\F(t)$ is convex and
\eqb\label{eq56b} \int\limits_{1}\limits^{\infty}\ H(t)\
\frac{dt}{t^2}\ =\ \infty\ . \eqe It was one of the most outstanding
results in the field of criteria for the solvability of the
degenerate Beltrami equations as it is clear from Theorem
\ref{th55}, see also Remark \ref{rmk555}.
\medskip

Subsequently, the fine theorems on the existence and uniqueness of
solutions in the Orlich--Sobolev classes have been established under
the condition (\ref{eq56b}) with the smooth $H$ and the condition
$tH'(t)\ge 5$, see Theorem 20.5.2 in the monograph \cite{AIM}, cf.
Lemma \ref{th555}, see also Remark  \ref{rmk6.333} above. However,
we have not found the work \cite{GMSV$_2$} in the reference list of
this monograph. The theorems on the existence and uniqueness of
solutions in the class $W^{1,2}_{loc}$ have been established  also
before it under $K_{\mu}(z)\le Q(z)\in W^{1,2}_{loc}$ in the work
\cite{MM}.
\medskip

Recently Brakalova and Jenkins have proved  the existence of ACL
ho\-meo\-mor\-phic solutions under (\ref{eq5.26V}), again with
$K_{\m}(z)$ instead of $K^T_{\m}(z,z_0)$, and with
\eqb\label{eq56aA} \F_{z_0}(t)\ \equiv\ \F(t)\ =\
h\left(\frac{t+1}{2}\right) \eqe where they assumed that $h$ is
increasing and convex and $h(x)\ge C_{\l}x^{\l}$ for any $\l > 1$
with some $C_{\l}>0$ and \eqb\label{eq5.26Lab}
\int\limits_{1}^{\infty}\frac{d\t}{\t h^{-1}(\t)}\ \ =\ \infty\ ,
\eqe see \cite{BJ$_2$}. Note that the conditions $h(x)\ge
C_{\l}x^{\l}$ for any $\l > 1$, in particular, under the above
sub-exponential integral constraints, see (\ref{eq555c}), imply that
$K_{\m}$ is locally integrable with any degree $p\in [1,\infty)$,
see (\ref{eq1.8}).
\medskip

Some of the given conditions are not necessary as it is clear from
the results in Section 4 and from the following lemma and remarks.

\blem{} \label{lem6.A} There exist continuous increasing convex
functions $\F :[1,\infty )\to [1,\infty )$ such that
\eqb\label{eq300A} \int\limits_{1}^{\infty} \log \F(t)\
\frac{dt}{t^2}\ =\ \infty\ , \eqe \eqb\label{eq300B}
\liminf\limits_{t\to\infty} \ \frac{\log \F(t)}{\log t}\ =\ 1\eqe
and, moreover, \eqb\label{eq300V} \F(t)\ge t\ \ \ \ \ \ \ \forall\
t\in [1,\infty )\ .\eqe \elem

Note that the examples from the proof of Lemma \ref{lem6.A} further
can be extended to $[0,\infty]$ by $\F(t)=t$ for $t\in[0,1]$ with
keeping all the given properties.

\brem\label{rmk6.12} The condition (\ref{eq300B}) implies, in
particular, that there exist no $\lambda
> 1, C_{\l} > 0$ and $T_{\l}\in[1,\infty )$ such that
\eqb\label{eq6.13} \Phi(t)\ \geq C_{\l}\, \cdot t^{\lambda}\ \ \ \ \
\ \forall\ t\geq T_{\l}\ .\eqe Thus, in view of Lemma \ref{lem6.A}
and Theorem \ref{th4.5A}, no of the conditions (\ref{eq6.13}) is
necessary in the existence theorems for the Beltrami equations with
the integral constraints of the type (\ref{eq5.27V}).

In addition, for the examples of $\F$ given in the proof of Lemma
\ref{lem6.A}, \eqb\label{eq6.14} \limsup\limits_{t\rightarrow
\infty}\ \frac{\log\,\Phi(t)}{\log\, t}\ =\ \infty\ , \eqe cf.
Proposition \ref{pr6.28} further. Finally, all the conditions
(\ref{eq333Y})--(\ref{eq333A}) from Theorem \ref{pr4.1aB} hold
simultaneously with (\ref{eq300A}) because the increasing convex
function $\F$ is absolutely continuous. \erem

{\it Proof of Lemma \ref{lem6.A}.} Further we use the known
criterion which says that a function $\Phi$ is convex on an open
interval $I$ if and only if $\Phi$ is continuous and its derivative
$\Phi^{\prime}$ exists and is non-decreasing in I except a countable
set of points in $I$, see e.g. Proposition $1.4.8$ in \cite{Bou}.We
construct $\Phi$ by induction sewing together pairs of functions of
the two types $\varphi(t)=\alpha+\beta t$ and $\psi(t)=ae^{bt}$ with
suitable positive parameters $a,b$ and $\beta$ and possibly negative
$\alpha.$ \medskip

More precisely, set $\Phi(t)=\varphi_1(t)$ for
$t\in[1,\gamma_1^{*}]$ and $\Phi(t)=\psi_1(t)$ for
$t\in[\gamma_1^{*},\gamma_1]$ where
$\varphi_1(t)=t,$\,\,$\gamma_1^{*}=e,$\,\,$\psi_1(t)=e^{-(e-1)}e^{t},$
\,\,$\gamma_1=e+1.$ Let us assume that we already constructed
$\Phi(t)$ on the segment $[1,\gamma_n]$ and hence that $\Phi(t)=a_n
e^{b_nt}$ on the last subsegment $[\gamma_n^{*},\gamma_n]$ of the
segment $[\gamma_{n-1},\gamma_n].$ Then we set
$\varphi_{n+1}(t)=\alpha_{n+1}+\beta_{n+1}t$ where the parameters
$\alpha_{n+1}$ and $\beta_{n+1}$ are found from the conditions
$\varphi_{n+1}(\gamma_{n})=\Phi(\gamma_n)$ and
$\varphi^{\prime}_{n+1}(\gamma_n)\geq \Phi^{\prime}(\gamma_n-0),$
i.e., $\alpha_{n+1}+\beta_{n+1}\gamma_{n}=a_ne^{b_n\gamma_n}$ and
$\beta_{n+1}\geq a_n b_n e^{b_n\gamma_n}.$ Let $\beta_{n+1}=a_n b_n
e^{b_n\gamma_n}$,
$\alpha_{n+1}=a_ne^{b_n\gamma_n}\left(1-b_n\gamma_n\right)$  and
choose a large enough $\gamma_{n+1}^{*}>\gamma_n$ from the condition
\eqb\label{eq6.15}
\log\,\left(\alpha_{n+1}+\beta_{n+1}\gamma_{n+1}^{*}\right)\,\leq\,\left(1+\frac1n\right)
log\,\gamma_{n+1}^{*} \eqe and, finally, set $\Phi(t)\equiv
\varphi_{n+1}(t)$ on $[\gamma_n,\gamma_{n+1}^{*}].$ \medskip

Next, we set $\psi_{n+1}(t)=a_{n+1}e^{b_{n+1}t}$ where parameters
$a_{n+1}$ and $b_{n+1}$ are found from the conditions that
$\psi_{n+1}(\gamma_{n+1}^{*})\ =\ \varphi_{n+1}(\gamma_{n+1}^{*})$
and $\psi_{n+1}^{\prime}(\gamma_{n+1}^{*})\ \ge\
\varphi_{n+1}^{\prime}(\gamma_{n+1}^{*})\ ,$ i.e.,
\eqb\label{eq6.16} b_{n+1}\ =\ \frac{1}{\gamma_{n+1}^{*}}\,\log\,
\frac{\alpha_{n+1}\ +\ \beta_{n+1}\gamma_{n+1}^{*}}{a_{n+1}} \eqe
and, taking into account (\ref{eq6.16}),
\eqb\label{eq6.17}b_{n+1}\ \geq\ \frac{\beta_{n+1}}{\alpha_{n+1}\
+\ \beta_{n+1}\gamma_{n+1}^{*}}\ . \eqe Note that (\ref{eq6.17})
holds if we take small enough $a_{n+1}>0$ in $(\ref{eq6.16})$. In
addition, we may choose here $b_{n+1}>1.$ \medskip

Now, let us choose a large enough $\gamma_{n+1}$ with
$e^{-1}\gamma_{n+1}\ \geq\ \gamma_{n+1}^{*}$ from the condition
that \eqb \log\,\psi_{n+1}\left(e^{-1}\gamma_{n+1}\right)\ \geq\
e^{-1}\gamma_{n+1}\ , \eqe i.e., \eqb\label{eq6.18} \log\,a_{n+1}\
+\ b_{n+1}e^{-1}\gamma_{n+1}\ \geq\ e^{-1}\gamma_{n+1}\ .\eqe Note
that (\ref{eq6.18}) holds for all large enough $\gamma_{n+1}$
because $b_{n+1}>1$ although $\log a_{n+1}$ can be negative.
\medskip

Setting $\Phi(t)=\psi_{n+1}(t)$ on the segment
$[\gamma_{n+1}^{*},\,\gamma_{n+1}],$ we have that
\eqb\label{eq6.19} \log\,\Phi(t)\ \geq\ t\ \ \ \ \ \ \forall
\,t\,\in\,[e^{-1}\gamma_{n+1},\gamma_{n+1}]\eqe where the
subsegment
$[e^{-1}\gamma_{n+1},\gamma_{n+1}]\subseteq[\gamma_{n+1}^{*},\gamma_{n+1}]$
has the logarithmic length 1. \medskip

Thus, (\ref{eq300V}) holds because by the construction $\F(t)$ is
absolutely continuous, $\F(1)=1$ and $\F'(t)\ge 1$ for all
$t\in[1,\infty);$ the equality (\ref{eq300A}) holds by
(\ref{eq6.19}); (\ref{eq300B}) by (\ref{eq300V}) and
(\ref{eq6.15}); (\ref{eq6.14}) by (\ref{eq6.19}).

\brem\label{rmk6.12V} Taking in the above construction in Lemma
\ref{lem6.A} $\b_{n+1}=1$ for all $n=1,2,\ldots $, $\a_{n+1}=b_n
e^{b_n\gamma_n}-\g_n$ and arbitrary $\g^*_{n+1}>\g_{n+1}$ we obtain
examples of absolutely continuous increasing functions $\F$ which
are not convex but satisfy (\ref{eq300A}), as well as all the
conditions (\ref{eq333Y})--(\ref{eq333A}) from Proposition
\ref{pr4.1aB}, and (\ref{eq300V}).

The corresponding examples of non-decreasing functions $\Phi$ which
are neither continuous, nor strictly monotone and nor convex in any
neighborhood of $\infty$ but satisfy (\ref{eq300A}), as well as
(\ref{eq333Y})--(\ref{eq333A}), and (\ref{eq300V}) are obtained in
the above construction if we take $\beta_{n+1}=0$ and
$\alpha_{n+1}>\gamma_n$ such that $\alpha_{n+1}
> \Phi(\gamma_n)$ and $\Phi(t)=\alpha_{n+1}$ for all $t\in
(\gamma_n,\gamma_{n+1}^{*}],$ $\gamma_{n+1}^{*}=\alpha_{n+1}.$ \erem

\bpr \label{pr6.28} Let $\Phi:[1,\infty)\rightarrow [1,\infty)$ be a
locally integrable function such that \eqb\label{eq6.29}
\int\limits_{1}^{\infty}\log{\Phi(t)}\ \frac{dt}{t^2}\ =\ \infty\ .
\eqe Then \eqb\label{eq6.30}
\limsup\limits_{t\rightarrow\infty}\frac{\Phi(t)}{t^{\lambda}}\,=\,\infty\
\ \ \ \ \ \forall\,\, \lambda\,\in\,{\Bbb R}\ .\eqe \epr

\brem\label{rmk6.31} In particular, (\ref{eq6.30}) itself implies
the relation (\ref{eq6.14}). Indeed, we have
from (\ref{eq6.30}) that there exists a monotone sequence
$t_n\rightarrow\infty$ as $n\rightarrow \infty$ such that
\eqb\label{eq6.32}\Phi(t_n)\,\geq\,t_n^n\,,\,\,\,\,n=1,2,\ldots,\eqe
i.e., \eqb\label{eq6.33} \frac{\log{\Phi(t_n)}}{\log{t_n}}\ \geq\
n,\,\,\,n=1,2,\ldots\,. \eqe \erem

{\it Proof of Proposition \ref{pr6.28}.} It is sufficient to
consider the case $\lambda > 0.$ Set $H(t)\,=\,\log{\Phi(t)},$ i.e.,
$\Phi(t)=e^{H(t)}.$ Note that $e^x\geq x^n/n!$ for all $x\geq 0$ and
$n=1,2\ldots\,,$ because $e^x=\sum\limits_{n=0}^{\infty}x^n/n!\,\,.$
Fix $\lambda>0$ and $n>\lambda.$ Then $q\colon =\lambda/n$ belongs
to $\left(0,1\right)$ and
$$\frac{H(t)}{t^q}\ \leq\ \left({\frac{\Phi(t)}{t^{\lambda}}}\right)^{\frac 1n}\cdot
\sqrt[n]{n!}\ .$$ Let us assume that \eqb\label{eq6.34} C\colon =\
\limsup\limits_{t\rightarrow\infty}\ \frac{\Phi(t)}{t^{\lambda}}\ <\
\infty\ .\eqe Then
$$\int\limits_{\triangle}^{\infty}H(t)\ \frac{dt}{t^2}\ <\ 2\sqrt[n]{Cn!}
\int\limits_{\triangle}^{\infty}\frac{dt}{t^{2-q}}\ =\ -\
\frac{2}{1-q}\ \frac{\sqrt[n]{Cn!}}
{t^{1-q}}\mid_{\triangle}^{\infty}\ =$$\,\,$$=\ \frac{2}{1-q}\
\frac{\sqrt[n]{Cn!}} {{\triangle}^{1-q}}\ <\ \infty$$ for large
enough $\triangle>1>0\,.$ The latter contradicts (\ref{eq6.29}).
Hence the assumption (\ref{eq6.34}) was not true and, thus,
(\ref{eq6.30}) holds for all $\lambda \in {\Bbb R}\,.$

\brem\label{rmk6.35} Lemma \ref{lem6.A} shows that, generally
speaking, $\limsup$ in (\ref{eq6.30}) cannot be replaced by $\lim$
for an arbitrary $\l>1$ under the condition (\ref{eq6.29}) even if
$\Phi$ is continuous, increasing and convex.  \erem

\bigskip

Applications of strong ring solutions to the theory of boundary
problems for the Beltrami equations will be published elsewhere, see
e.g. \cite{Dy}, cf. \cite{RS} and \cite{Lo}.

\medskip
\noindent {\bf Vladimir Ryazanov:}\\ Institute of Applied
Mathematics\\ and Mechanics, NAS of Ukraine, \\ 74 Roze Luxemburg
str.,\\ 83114, Donetsk, UKRAINE\\
Email: ${\tt vlryazanov1@rambler.ru}$
\medskip

\noindent {\bf Uri Srebro:}\\ Technion, \\ Haifa 32000, ISRAEL\\
Email: {\tt srebro@math.technion.ac.il}\\

\medskip
\noindent {\bf Eduard Yakubov:}\\ Holon Institute of Technology,\\
52 Golomb St., P.O.Box 305,\\ Holon 58102,
ISRAEL\\ Email: {\tt yakubov@hit.ac.il}\\

\end{document}